\documentclass[12pt,letter]{article}
\usepackage{amsmath,amsthm,enumerate}
\title{An extended class of minimax generalized Bayes estimators of 
regression coefficients}
\author{Yuzo Maruyama\thanks{University of Tokyo \quad e-mail:
maruyama@csis.u-tokyo.ac.jp} 
\and William, E. Strawderman \thanks{Rutgers University \quad e-mail: straw@stat.rutgers.edu}} 

\newtheorem{thm}{Theorem}[section]
\newtheorem{lemma}{Lemma}[section]
\newtheorem{corollary}{Corollary}[section]

\allowdisplaybreaks
\date{}

\begin{document}
\maketitle
\begin{abstract}
We derive minimax generalized Bayes estimators of regression
 coefficients in the general linear model with spherically symmetric
 errors under invariant quadratic loss for the case of unknown scale.
The class of estimators generalizes the class considered in Maruyama and
 Strawderman (2005) to include non-monotone shrinkage functions. \\

\end{abstract}

\section{Introduction}
In this paper we consider minimax generalized Bayes estimators of 
the regression coefficients in 
the general linear model with homogeneous spherically symmetric errors.
We start with the familiar linear regression model
$ Y=A\beta+\epsilon $ where $ Y$ is an $ N \times 1$
vector of observations, $A$ is the known $N \times p$ design matrix of
rank $p$, $ \beta$ is the $p \times 1$ vector of unknown regression 
coefficients, and $ \epsilon $ is an $N \times 1$ vector of experimental
errors. We assume $ \epsilon $ has a spherically symmetric distribution
with a density $ \sigma^{-N}f(\epsilon'\epsilon/\sigma^2)$,
where $ \sigma $ is an unknown scale parameter and $ f(\cdot)$ is a
nonnegative function on the nonnegative real line, which satisfies
\begin{align*}
 \int_0^\infty t^{N/2-1}f(t)dt <\infty .
\end{align*}
The problem is to estimate $\beta$.
The least squares estimator of $ \beta$ is $ \hat{\beta}=(A'A)^{-1}A'y$.
In order to treat the estimation problem 
from the decision-theoretic point of view,
we measure the loss in estimating $ \beta$ by $b$ with so called
scale invariant ``predictive loss'' functions
\begin{equation}
L(b,\beta, \sigma^2)=\sigma^{-2}(b-\beta)'A'A(b-\beta). \label{gloss}
\end{equation}
Then the risk function of an estimator $b$ is given by
$ R(b,\beta,\sigma^2)=E[L(b,\beta,\sigma^2)]$.
The least squares estimator $ \hat{\beta}$ is minimax with constant risk.
Therefore, $ b$ is a minimax estimator of $ \beta$ if and only if
$ R(b,\beta,\sigma^2) \leq R(\hat{\beta},\beta, \sigma^2) $ for all 
$\beta$ and $ \sigma^2$,
the search for estimators better than $ \hat{\beta}$ is a search
for minimax estimators.
It is, of course, well known that estimators dominating $\beta$ exist
when $p \geq 3$.
In this paper, we study generalized Bayes minimax estimators.

To simplify expressions and to make matters a bit clearer 
it is helpful
to rotate the problem via the following transformation,
so that the covariance matrix of $ \hat{\beta}$, $\sigma^2(A'A)^{-1}$ 
becomes diagonal.
Let $P$ be the $p \times p$ orthogonal matrix of 
eigenvectors of $(A'A)^{-1}$, 
with  $ d_1 \geq d_2 \geq \dots \geq d_p$ as eigenvalues,
it follows that
\begin{equation*}
 P'(A'A)^{-1}P=D, \quad P'P=I_p
\end{equation*}
where $D=\mbox{diag}(d_1, \dots ,d_p )$.
Let also $Q$ be an $N \times N$ orthogonal matrix
such that
\begin{equation*}
QA=\begin{pmatrix}
          D^{-1/2}P' \\
          0
         \end{pmatrix}.
\end{equation*}
Next, define two random vectors
$ X=(X_{1},\cdots,X_{p})'$ and $ Z=(Z_{1},\cdots,Z_{n})' $
where $n=N-p$ by
\begin{equation*}
\begin{pmatrix}
X\\
Z
\end{pmatrix}
=QY. 
\end{equation*}
Then $(X',Z')'$ has the joint density given by
\begin{equation}
\label{f}
\sigma^{-p-n}f(\{(x -\theta)'(x-\theta)+z'z\}
/\sigma^{2}),
\end{equation}
where $ \theta=P'\beta$.
Notice also that
$ X $ and $ Z'Z$ can be expressed  as $ D^{-1/2}P'\hat{\beta} $ and
$(y-A\hat{\beta})'(y-A \hat{\beta})$, respectively.
Also, as is customary, denote $Z'Z$ by $S$.
We assume throughout this paper that $p \geq 3$ and $n
\geq 1$.

The original problem is thus equivalent to estimation of $ \theta $ under the
loss function 
$  L(\delta,\theta, \sigma^2)=(\delta - \theta)'(\delta -
\theta) /\sigma^2$. 
We will consider the problem in this equivalent canonical form.

This paper is best viewed as a companion paper to
Maruyama and Strawderman (2005). In that paper a class of minimax
generalized Bayes 
estimators was derived for the above canonical problem 
when the errors were normally distributed, and a subclass was shown to
be generalized Bayes and minimax for the entire class of spherically
symmetric error distributions.

In this paper we enlarge the class of generalized Bayes minimax
estimators for each of above classes of distributions. 
To do so we must first enlarge the class of
minimax estimators to accommodate non-monotone shrinkage functions.
Section 2 is denoted to this generalization. Section 3 develops an
extended class of generalized Bayes minimax estimators for spherical
normal error distributions. To a large degree, the results therein
represent extentions to the case of unknown variance, of the results of
Maruyama (1998, 2004) in the known variance case.
Section 4 extends these results to the case of general spherically symmetric
error distributions with unknown scale. A subclass of the estimators
in Section 3 is shown to be generalized Bayes and minimax for the entire
class of spherical error distributions simultaneously (subject to
finiteness of moments).

\section{Minimaxity}
In this section, we give a sufficient condition for minimaxity in the
general spherically symmetric case. 
\begin{thm}
Suppose $(X',Z')'$ has a distribution given by \eqref{f}. 
Then $ \delta_{\phi}$ given by
\begin{align} \label{shrinkage}
 \delta_\phi=\left\{1-\frac{S}{\|X\|^2}
\phi\left(\frac{\|X\|^2}{S}\right)\right\}X
\end{align}
is minimax 
if 
\begin{align} \label{sufficient-minimax}
 \frac{\phi(w)}{w}\left\{
(n+2)\phi(w)-2(p-2)
\right\}-4\phi'(w)\{1+\phi(w)\} \leq 0.
\end{align}
\end{thm}
The interesting point of the theorem is that the sufficient condition
for minimaxity does not depend on $f$.
Such distributional robustness has already been noted in the literature,
but only the following tractable subset of the above result is typically used:
\begin{corollary}\label{cor-minimax-naive}
 If $\phi$ is monotone nondecreasing and $0\leq \phi(w) \leq 2(p-2)/(n+2)$,
$\delta_\phi$ is minimax in the general
 spherically symmetric case.
\end{corollary}
For completeness we give the proof.
We use the version derived in Kubokawa and Srivastava (2001)
but earlier versions appear in Robert (1994) and elsewhere.
\begin{proof}
Let 
\begin{equation*}
 F(x)=\frac{1}{2}\int_{x}^{\infty}f(t)dt
\end{equation*}
and define
\begin{align*}
 E^{f}
[h(X,Z)]&=\int \int h(x,z)\sigma^{-N}
f \left( \frac{(x - \theta)'(x - \theta)}{\sigma^2}+
\frac{z'z}{\sigma^2}\right)dx dz \\
 E^{F}
[h(X,Z)]&=\int \int h(x,z)\sigma^{-N}
F \left( \frac{(x - \theta)'(x - \theta)}{\sigma^2}+
\frac{z'z}{\sigma^2}\right)dx dz,
\end{align*}
where $h(x,z)$ is an integrable function. 
Note that
$F(x)$ is not necessarily a probability density, and hence
$E^F$ is not necessarily an expectation symbol in any strict sense.

The identities corresponding to the Stein 
and chi-square identities for the normal distribution (Stein (1973) and
Efron and Morris (1976)),
\begin{align}
  E^{f}[(X_i - \theta_i )h(X,Z)]&=  \sigma^2 
E^{F}[(\partial/\partial X_i)h(X,Z)], \label{steiniden} \\ 
E^{f}[S g(S)]&=\sigma^2 E^{F}[n g(S)+2S
 g'(S)],\label{chiiden} 
\end{align}
where $S=Z'Z$, are useful in our proof.

The risk of $ \hat{\theta}_{\phi}$ is given by
\begin{align}
 R(\theta, \sigma^2, \hat{\theta}_{\phi}) 
&= E^f \left[ ( \hat{\theta}_{\phi} - \theta)'( \hat{\theta}_{\phi} -
 \theta) /\sigma^2 \right] \nonumber \\
&=R(\theta,\sigma^2,X)
 +E^f \left[ \frac{S^2}{\sigma^2 \|X\|^2}
\phi^2 \left(  \frac{\|X\|^2}{S} \right)\right] \nonumber \\
&\quad -2E^f \left[ \sum \frac{S}{\sigma^2 \|X\|^2}X_i(X_i - \theta_i)
 \phi\left( \frac{\|X\|^2}{S} \right) \right]. \label{risk1}
\end{align}
Let $W=\|X\|^2/S$.
For the second term in \eqref{risk1}, using \eqref{chiiden},
we have
\begin{align*}
E^f \left[ 
\frac{S}{\sigma^2\|X\|^2} \left\{ S 
\phi^2 \left(  \frac{
\|X\|^2}
{S} \right)\right\}\right] 
= E^F \left[  
 (n+2)\frac{\phi^{2}(W)}{W}
-4
\phi(W)\phi'(W)  \right].
\end{align*}
For the third term in \eqref{risk1}, using \eqref{steiniden}, we have
\begin{align*}
&\sum E^f \left[ \frac{1}{\sigma^2}(X_i - \theta_{i})X_i
\left( \frac{\|X\|^2}{S}\right)^{-1}\phi\left( \frac{\|X\|^2}{S} \right) \right] \\
&=E^F \left[ p\frac{\phi(W)}{W}+2\frac{
\|X\|^2}{S}
\left\{ \frac{\phi'(W)}{W}-\frac{\phi(W)}{W^2}\right\} \right].
\end{align*}
Hence  we have 
\begin{align*}
& R(\theta,\sigma^2,\hat{\theta}_{\phi})- R(\theta,\sigma^2,X)\\
& =
E^F \left[
\frac{\phi(W)}{W}\left\{
(n+2)\phi(W)-2(p-2)
\right\}-4\phi'(W)(1+\phi(W)) \right],
\end{align*}
which completes the proof.
\end{proof}
An immediate corollary, which we believe is often more tractable,
is given by dividing the left hand side of the inequality
\eqref{sufficient-minimax} by $\{1+\phi(w)\}\phi(w)/w$.
\begin{corollary} \label{cor-minimax}
\begin{enumerate}[(i)]
\item  $\delta_\phi$ given by \eqref{shrinkage}
is minimax if
\begin{align} \label{tractable-c}
\frac{(n+2)\phi(w)-2(p-2)}{1+\phi(w)}-4w\frac{\phi'(w)}{\phi(w)} \leq 0.
\end{align}
\item
Suppose $0 \leq \phi(w)\leq M_1$ and $ w\phi'(w)/\phi(w)\geq -M_2$.
The estimator $\delta_\phi$ given by \eqref{shrinkage}
is minimax if
\begin{align*}
\frac{(n+2)M_1-2(p-2)}{1+M_1}+4M_2 \leq 0.
\end{align*}
\end{enumerate}
\end{corollary}
It is clear that part (ii) of the corollary above allows $\phi$ be
non-monotonic. 
Actually we show in the next section that
there exists a class of minimax generalized Bayes estimators with
non-monotone $\phi$.
\section{An extended class of Generalized Bayes minimax estimators in the normal case}
In this section, we extend a class of generalized Bayes minimax
estimators in Maruyama and Strawderman (2005) in the normal case.

Suppose the sampling distribution of $(X',Z')'$ is normal with covariance
matrix $ \sigma^2 I_n  $
and mean vector $(\theta',0')'$.
As in Maruyama and Strawderman (2005), 
the class of hierarchical priors we consider is as follows;
\begin{align} \label{prior}
 \theta|\eta,\lambda \sim N_p(0,\eta^{-1}(1-\lambda)\lambda^{-1} I_p), 
\ \eta \sim \eta^e, \lambda \sim \lambda^a(1-\lambda)^{b} .
\end{align}
Then the marginal density of $X$, $S$, $ \lambda$ and $ \eta$ is proportional to
\begin{align}
& \int_{R^p}\exp\left( -\frac{\eta}{2} \left\{ \|x-\theta\|^2
+\frac{\lambda}{1-\lambda}\|\theta\|^2
\right\} -\frac{\eta s}{2}
\right)\frac{\eta^{p/2+n/2+e}\lambda^{p/2+a}}{(1-\lambda)^{p/2-b}}
d \theta \nonumber \\
& \propto \exp \left( -\frac{\eta s}{2}
(1+\lambda w) \right)
\eta^{p/2+n/2+e}\lambda^{p/2+a}  
(1- \lambda)^{b}, \label{marginal}
\end{align}
where $ w=\|x\|^2/s $.
Under quadratic loss,
the generalized Bayes estimator for such a hierarchical priors can be
expressed as 
\eqref{marginal},  
\begin{equation} \label{estimator-GB}
\hat{\theta}=\frac{E(\eta \theta|X,S)}{E(\eta|X,S)} =
\left(1 - \frac{E(\lambda
	 \eta|X,S)}{E(\eta|X,S)}\right)X
=\left( 1 -\frac{\phi(W)}{W}\right)X.
\end{equation}
When $  p/2+n/2+e+2>0$,
\begin{equation}
 \int_0^{\infty} \eta^{p/2+n/2+e+1}\exp \left( -\frac{\eta s}{2}
(1+\lambda w) \right)d\eta \propto
(1+\lambda w)^{-p/2-n/2-e-2}, 
\end{equation}
hence
\begin{align}
\phi(w)
=w\frac{\int_0^{1}\lambda^{p/2+a+1}(1-\lambda)^{b}(1+w
\lambda)^{-p/2-n/2-e-2}dt}
{\int_0^{1}\lambda^{p/2+a}(1-\lambda)^{b}(1+w\lambda )^{-p/2-n/2-e-2}dt} , \label{phi2}
\end{align}
which is well-defined for $ a > -p/2-1$ and $ b> -1$. 

Here is a result of Maruyama and Strawderman (2005).
\begin{thm}[Maruyama and Strawderman (2005)]
Suppose $ b \geq 0$, $ e> -p/2-n/2-1$ and $ -p/2-1< a <n/2+e $.
Then $\phi(w)$ is monotone increasing and approaches  $
 (p/2+a+1)/(n/2+e-a)$  as $w \to \infty$.
Hence
$ \delta_\phi$ with $\phi$ given by \eqref{phi2}
is minimax if
\begin{equation*}
 0 \leq \frac{p/2+a+1}{n/2+e-a} \leq 2\frac{p-2}{n+2}.
\end{equation*} 
\end{thm}
The key assumption of the theorem is $b \geq 0$ while 
the estimator with $\phi(w)$ itself is well-defined when $b>-1$.
In fact, there is a big difference between the two cases: $ b \geq 0$ and
$-1< b <0$. An immediately apparent difference is 
that $(1-\lambda)^{b}$  is unbounded for $ -1 <b<0$
while
it is bounded for $b \geq 0$. Technically, because of this unboundedness,
the standard integration by parts technique
fails to work well when $-1 <b<0$.
Furthermore the shrinkage factor $\phi(w)$ for $-1<b<0$ is not often monotonic,
which means the tractable sufficient condition for minimaxity given by
Corollary \ref{cor-minimax-naive} is not applicable.

We will see below that 
proving minimaxity with $b<0$ 
requires a different approach from that in
Maruyama and Strawderman (2005). 
Actually it will be done through the expression of $\phi$
by 
hypergeometric functions 
\begin{align}
F(\alpha,\beta;\gamma;z)=1+\sum_{i=1}^\infty\frac{(\alpha)_i
 (\beta)_i}{(\gamma)_i}\frac{x^i}{i!} \label{HGF}
\end{align}
where $(\alpha)_i=\alpha \cdots (\alpha+i-1)$.
The following lemma summarizes the relationships we use.
All formulas are from Abramowitz and Stegun (1964)
and the number following AS in each expression below is the formula
number in Abramowitz and Stegun (1964).  
\begin{lemma} \label{abra}
\begin{itemize}
\item AS.15.3.1
\begin{align}\label{3.1}
F(\alpha,\beta;\gamma;z)=\frac{\Gamma(\gamma)}{\Gamma(\beta)\Gamma(\gamma-\beta)}
\int_0^1 t^{\beta-1}(1-t)^{\gamma-\beta-1}(1-tz)^{-\alpha}dt 
\end{align}
when $\gamma>\beta>0$.
\item AS.15.3.4
\begin{align} \label{3.4}
 F(\alpha,\beta;\gamma;z)=(1-z)^{-\alpha}F(\alpha,\gamma-\beta;\gamma;z/(z-1)) .
\end{align}
\item AS.15.2.25
\begin{align}\label{2.25}
&\gamma(1-z)F(\alpha, \beta ; \gamma ;z)-\gamma F(\alpha,\beta-1;\gamma;z)
\\
& \qquad \qquad +z(\gamma-\alpha)F(\alpha,\beta;\gamma+1;z)=0. \notag
\end{align}
\item 
AS.15.2.18
\begin{align}\label{2.18}
&(\gamma-\alpha-\beta)F(\alpha,\beta;\gamma;z)-(\gamma-\alpha)
F(\alpha-1,\beta;\gamma;z)\\
& \qquad \qquad +\beta(1-z)F(\alpha,\beta+1;\gamma;z)=0. \notag
\end{align}
\item AS.15.1.20
\begin{align} \label{1.20}
\lim_{z \to 1} F(\alpha,\beta;\gamma;z)
=\frac{\Gamma(\gamma)\Gamma(\gamma-\alpha-\beta)}
{\Gamma(\gamma-\alpha)\Gamma(\gamma-\beta)}
\end{align}
when $\gamma \neq 0,-1,-2, \dots $ and $\gamma-\alpha-\beta>0$.
\item AS.15.3.10 
\begin{align}\label{3.10}
\lim_{z \to 1}\{ -\log(1-z)\}^{-1}
 F(\alpha,\beta;\alpha+\beta;z)=\frac{\Gamma(\alpha+\beta)}{\Gamma(\alpha)\Gamma(\beta)}.
\end{align}
\item AS.15.3.3
\begin{align} \label{3.3}
 F(\alpha,\beta;\gamma;z)=(1-z)^{\gamma-\alpha-\beta}
F(\gamma-\alpha,\gamma-\beta;\gamma;z).
\end{align}
\item AS.15.2.1
\begin{align} \label{2.1}
 (d/dz)F(\alpha, \beta ; \gamma ; z)=\{\alpha \beta/\gamma \}
F(\alpha+1,\beta+1 ; \gamma+1;z).
\end{align}
\end{itemize}
\end{lemma}

We first use Lemma \ref{abra} to re-express $\phi(w)$.
\begin{lemma} \label{phi-G}
Let $\phi$ be given by \eqref{phi2}. Then provided $a>-p/2-1$, $b>-1$ and $
p/2+n/2+e+2>0$, $\phi(w)$ is expressed as
\begin{align}
\phi(w) 
= \frac{1-G(v)}{\frac{n/2+e-a}{p/2+a+1}+G(v)} \label{phi--2}
\end{align}
where 
\begin{align*}
 G(v)=\frac{F(b,p/2+n/2+e+1;p/2+a+b+2;v)}{F(b+1,p/2+n/2+e+1;p/2+a+b+2;v)}. \label{Gv}
\end{align*}
\end{lemma}
\begin{proof}
Using \eqref{3.1} and \eqref{3.4},
we have 
\begin{align}
& \int_0^1 t^{\beta-1}(1-t)^{\gamma-\beta-1}(1+tz)^{-\alpha}dt \\
& \quad =
\frac{\Gamma(\beta)\Gamma(\gamma-\beta)}{\Gamma(\gamma)}
(1+z)^{-\alpha}F(\alpha,\gamma-\beta;\gamma;z/(1+z)) \notag
\end{align}
when $\gamma>\beta>0$ and hence by \eqref{phi2}
\begin{align}
 \phi(w)=\frac{v}{1-v}\frac{p/2+a+1}{d}
\frac{F(b+1,c;d+1;v)}{F(b+1,c;d;v)}, \label{phi--1}
\end{align}
where $v=w/(w+1)$, $c=p/2+n/2+e+2$, and $d=p/2+a+b+2$.
Using \eqref{2.25} and \eqref{2.18},
we have
\begin{align}
\phi(w) &=
 -1+\frac{F(b+1,c-1;d;v)}{(1-v)F(b+1,c;d;v)} \notag \\
&=-1+\frac{c-1}{n/2+e-a+(p/2+a+1)G(v)} \notag \\
&= \frac{1-G(v)}{\frac{n/2+e-a}{p/2+a+1}+G(v)} .
\end{align}
\end{proof}
The next lemma gives properties of ratios of hypergeometric functions
such as $G$, which we employ in demonstrating minimaxity.

\begin{lemma}\label{phi-K}
Let 
\begin{align}
K(v)=\frac{F(b,\beta;\gamma;v)}{F(b+1,\beta;\gamma;v)}  \label{KK}
\end{align}
for $\beta>0$ and $ 0<b+1<\gamma < b+1+\beta $.
Then,
\begin{enumerate} [(i)]
\item  $K(0)=1$
\item 
$\lim_{v \to 1}K(v)=0$.
\item If $b\geq 0$ or if $ -1 < b<0$ and $\gamma \geq \beta $, 
$K(v)$ is monotone decreasing.
\item If  $ -1 <b<0$ and $ \gamma < \beta $,
the minimum of $K(v)$ is a negative value and $K(v)$
approaches $0$ from the below as $v$ approaches $1$.
Also $\inf_{0\leq v \leq 1}K(v) \geq b/(b+1)$.
\end{enumerate}
\end{lemma}
\begin{proof}
First note that $K(0)=1$ by \eqref{KK}. 
Next by \eqref{1.20} and \eqref{3.3},
\begin{align} \label{1.20-3.3}
\lim_{z \to 1}(1-z)^{\alpha+\beta-\gamma} F(\alpha,\beta;\gamma;z)=
 \frac{\Gamma(\gamma)\Gamma(\alpha+\beta-\gamma)}{\Gamma(\alpha)\Gamma(\beta)} 
\end{align}
when $\gamma >0 $ and $\gamma-\alpha-\beta<0$. Hence we have by 
\eqref{1.20-3.3}, \eqref{3.10} and \eqref{1.20} respectively,
\begin{align} 
\label{F-approx-1} F(b,\beta;\gamma;v) &\approx
\begin{cases}
 (1-v)^{\gamma-b-\beta}\frac{\Gamma(\gamma)\Gamma(b+\beta-\gamma)}
{\Gamma(b)\Gamma(\beta)} & \mbox{ if } \gamma-b-\beta<0 \\
\{-\log(1-v)\}\frac{\Gamma(\gamma)}{\Gamma(b)\gamma(\beta)} & \mbox{ if }\gamma-b-\beta=0 \\
\frac{\Gamma(\gamma)\Gamma(\gamma-b-\beta)}
{\Gamma(\gamma-b)\Gamma(\gamma-\beta)} & \mbox{ if } 0<\gamma-b-\beta<1 \\
\end{cases}
\end{align}
and also by \eqref{1.20-3.3}
\begin{align}
 F(b+1,\beta;\gamma;v) &\approx
 (1-v)^{\gamma-b-1-\beta}\frac{\Gamma(\gamma)\Gamma(b+1+\beta-\gamma)}
{\Gamma(b+1)\Gamma(\beta)} \mbox{ if } \gamma-b-\beta<1 \label{F-approx-2}
\end{align}
where $f(v) \approx g(v)$ means $\lim_{v \to 1} f(v)/g(v)=1$.
We easily see that the ratio $  F(b,\beta;\gamma;v)/ F(b+1,\beta;\gamma;v)=K(v)$
goes to zero as $v \to 1$ under the assumption $\gamma-b-\beta <1$.
Hence part (ii) follows.

When $b\geq 0$, 
$K(v)$ is decreasing from the monotone likelihood ratio property
of the kernel of $k(v)$. Hence the first assertion of part (iii) follows.

When $ -1 < b<0$ and $ \gamma - \beta \geq 0$, the numerator $ F(b,\beta;\gamma;v)$
is always positive because it can be rewritten as
$ (1-v)^{\gamma-b-\beta}F(\gamma-b,\gamma-\beta;\gamma;v) $.
Also the numerator $ F(b,\beta;\gamma;v)$ with $-1<b<0$ is decreasing in $v$ and the
positive denominator $ F(b+1,\beta;\gamma;v)$ is increasing in $v$.
Hence the second assertion of part (iii) follows.

To show part (iv), note that $ \Gamma(x)>0$ if $x>0$ and $ \Gamma(y)<0$ if
 $-1<y<0$. By assumption 
$-1<b<0$, $\gamma < \beta$, $\beta>0$, 
$0 <b+1<\gamma <b+1+\beta$.
It then follows using the additional assumptions in the first and third
 lines of \eqref{F-approx-1},
that there is exactly 1 negative factor in each constant term and hence
each constant term is negative. Since, also, the denominator of $K(v)$,
$F(b+1,\beta;\gamma;v)$,
is positive, it follows that $K(v)$ approaches 0 from below as $v$
 approaches 1.
Let $K(v)$ take on its minimum value at $v_0$.
Using the formula
\begin{align*}
 (A/B)'= (B'/B)\{ A'/B'- A/B \}
\end{align*}
and \eqref{2.1}, we have
 $K(v_0)=K_1(v_0)$
where
\begin{align*}
K_1(v)=
\frac{b}{b+1}\frac{F(b+1,\beta+1;\gamma+1;v)}{F(b+2,\beta+1;\gamma+1;v)}. 
\end{align*}
Since $b/(b+1)<0$ and $b+1>0$, 
$K_1(v)$ is increasing in $v$.
Therefore
\begin{align}
K(v) \geq K(v_0)= K_1(v_0) >K_1(0)=b/(b+1). \label{lower-bound-1}
\end{align}
Hence part (iv) follows.
\end{proof}
The following corollary gives the behavior of $\phi(w)$ and follows
immediately from Lemma \ref{phi-G} and \ref{phi-K}.
\begin{corollary} \label{upper-bound-phi}
Assume that $ e>-p/2-n/2-1$, $-p/2-1 <a<n/2+e $ and $b>-1$. Then 
\begin{enumerate}[(i)]
\item $  \lim_{w \to \infty} \phi(w)=\alpha$
where $ \alpha=(p/2+a+1)/(n/2+e-a)$.
\item When $b \geq \min(n/2+e-a-1, 0)$, $ \phi(w)$ is monotone increasing.
\item When $ -1<b<\min(n/2+e-a-1, 0)$, $\phi(w)$ is not monotonic.
\item When 
\begin{gather} \label{restriction-b-1}
-\frac{n/2+e-a}{p/2+n/2+e+1}<b<\min(n/2+e-a-1, 0), \\
0 \leq \phi(w) \leq
\frac{1-b/(b+1)}{1/\alpha+b/(b+1)}
=\frac{(p/2+a+1)}{n/2+e-a+b(p/2+n/2+e+1)}. \notag
\end{gather}
\end{enumerate}
\end{corollary}
Note: By analogy with the known variance case (Alam (1973)
and Maruyama (2003, 2004)), 
it may be expected that the choice $-1<b<0$ would lead to a non-monotone
$\phi$. 
However part (ii) of the corollary shows that this need not be true in the
unknown variance case,
and in fact monotonicity depends
on the relationship between $n$, $a$ and $e$.
The addition of the restriction \eqref{restriction-b-1} 
in part (iv) is necessitated by the fact
that
$ \lim_{b \to -1} b/(b+1) = -\infty $. Hence a value of $b$ close to $-1$ would
cause
the upper bound on $\phi$, $\{1-b/(b+1)\}/\{1/\alpha+b/(b+1) \}$ to be negative.
Thanks to the restriction, the denominator is positive.

The next result gives a lower bound for $w\phi'(w)/\phi(w)$.
\begin{lemma} \label{lower-bound-phi}
\begin{align}
w\frac{\phi'(w)}{\phi(w)} \geq \frac{(p/2+a+2)b}{b+1}
\end{align}
provided that $e>-p/2-n/2-1$, $ -p/2-1<a<n/2+e-a$ and $ -1 < b< \min(n/2+e-a-1,0)$.
\end{lemma}
\begin{proof}
Note 
\begin{align*}
\{wA(w)/B(w)\}' /\{A(w)/B(w)\}= 1+wA'(w)/A(w)-wB'(w)/B(w).
\end{align*}
Hence
\begin{align*}
& \frac{w\phi'(w)}{\phi(w)}=1+c w \left\{ -
\frac{\int_0^1 \lambda^{p/2+a+2}(1-\lambda)^b
(1+w\lambda)^{-c-1}d\lambda}
{\int_0^1 \lambda^{p/2+a+1}(1-\lambda)^b
(1+w\lambda)^{-c}d\lambda} \right.\\
& \qquad \qquad \quad +
\left. \frac{\int_0^1 \lambda^{p/2+a+1}(1-\lambda)^b
(1+w\lambda)^{-c-1}d\lambda}
{\int_0^1 \lambda^{p/2+a}(1-\lambda)^b
(1+w\lambda)^{-c}d\lambda}
\right\}
\end{align*}
where $c=p/2+n/2+e+2$. 
Using \eqref{3.1} and \eqref{3.4},
we have
\begin{align*}
\frac{w\phi'(w)}{\phi(w)} 
& =-c v \left\{ 
\frac{\frac{p}{2}+a+2}{d+1}
\frac{F(b+1,c+1;d+2;v)}
{F(b+1,c;d+1;v)} \right. \\
& \quad - \left. \frac{\frac{p}{2}+a+1}{d}\frac{F(b+1,c+1;d+1;v)}
{F(b+1,c;d;v)} \right\}+1,
\end{align*}
where $ v=w/(w+1)$ and $d=p/2+a+b+2$.
Using \eqref{2.25} and \eqref{2.18},
\begin{align*}
 \frac{w\phi'(w)}{\phi(w)} 
& =c(1-v) \left\{ 
\frac{F(b+1,c+1;d+1;v)}
{F(b+1,c;d+1;v)}
-  \frac{F(b+1,c+1;d;v)}
{F(b+1,c;d;v)} \right\}+1 \\
& = (\frac{p}{2}+a+2)
\frac{F(b,c;d+1;v)}
{F(b+1,c;d+1;v)} 
-  (\frac{p}{2}+a+1)\frac{F(b,c;d;v)}
{F(b+1,c;d;v)}.
\end{align*}
Note $ F(b,c;d+1;v) > F(b,c;d;v)$ since $-1<b<0$.
For $v$ which satisfies $ F(b,c;d;v)>0$, 
we have
\begin{align*}
 \frac{w\phi'(w)}{\phi(w)} 
=F(b,c;d;v) \left( \frac{p/2+a+2}{F(b+1,c;d+1;v)} -
\frac{p/2+a+1}{F(b+1,c;d;v)} \right) \geq 0,
\end{align*}
because $ F(b+1,c;d;v) > F(b+1,c;d+1;v)$.
For $v$ which satisfies 
$ F(b,c;d+1;v) > 0 > F(b,c;d;v)$, $ w\phi'(w)/\phi(w) $ is clearly
 nonnegative.
For $v$ which satisfies $F(b,c;d+1;v) < 0$ and $ -1 < b< \min(n/2+e-a-1,0)$,
we have 
\begin{align*}
 \frac{w\phi'(w)}{\phi(w)} \geq
(\frac{p}{2}+a+2)
\frac{F(b,c;d+1;v)}
{F(b+1,c;d+1;v)} \geq (\frac{p}{2}+a+2)\frac{b}{b+1}
\end{align*}
by Lemma \ref{phi-K} (iv).
\end{proof}
The main result is the following.
\begin{thm} \label{main-theorem}
Suppose $e >-p/2-n/2-1$. 
\begin{enumerate}[(i)]
\item \ [monotone $\phi$] When
\begin{align} \label{m-condition-a}
 -p/2-1<a \leq \frac{c(p,n)(n/2+e)-p/2-1}{1+c(p,n)}
\end{align}
where $c(p,n)=2(p-2)/(n+2)$ and $ b \geq \min(n/2+e-a-1,0)$,
the generalized Bayes estimator is minimax.
\item \ [non-monotone $\phi$] 
When
\begin{align} \label{m-condition-a-1}
 -p/2-1<a < \frac{c(p,n)(n/2+e)-p/2-1}{1+c(p,n)}
\end{align}
\begin{align} \label{main-b}
-1 <-(n+2)\frac{c(p,n)-\alpha}{4(p+a+1)(\alpha+1)} \leq b < \min(0,n/2+e-a-1)
\end{align}
where $\alpha=(p/2+a+1)/(n/2+e-a)=\lim_{w \to \infty}\phi(w)$,
 the generalized Bayes estimator 
is minimax.
\end{enumerate} 
\end{thm}
\begin{proof}
First we prove part (i). 
Monotonicity of $\phi$ follows from Corollary \ref{upper-bound-phi} (ii)
 since 
$b \geq \min(0, n/2+e-a-1)$. 
Also using Corollary \ref{upper-bound-phi} (i) and  \eqref{m-condition-a}
\begin{align*}
 0 \leq \phi \leq \alpha =\frac{p/2+a+1}{n/2+e-a} \leq 2\frac{p-2}{n+2}=c(p,n).
\end{align*}
Hence minimaxity follows from Corollary \ref{cor-minimax-naive} and part
(i) follows.

Next we consider part (ii), the non-monotonic $\phi$ case.
The lower bound in \eqref{main-b} 
\begin{align*}
 -(n+1)\frac{c(p,n)-\alpha}{4(p+a+1)(\alpha+1)} 
\end{align*}
is negative because of \eqref{m-condition-a-1}.
Also Corollary \ref{upper-bound-phi} (i)  and \eqref{m-condition-a-1}
implies
\begin{align*}
0< \lim_{w \to \infty}\phi(w)=\alpha =\frac{p/2+a+1}{n/2+e-a} < 2\frac{p-2}{n+2}=c(p,n).
\end{align*}
Also since the lower bound in \eqref{main-b} satisfies
\begin{align*}
& -(n+2)\frac{c(p,n)-\alpha}{4(p+a+1)(\alpha+1)}
-\left(-\frac{n/2+e-a}{p/2+n/2+e+1}\right) \\
& \geq \frac{1}{\alpha+1}\left( 1-
 \frac{(n+2)c(p,n)}{4(p+a+1)}\right) \\
& = \frac{1}{\alpha+1}\left( 1-\frac{p-2}{p+2(p/2+a+1)}\right) \\
& \geq \frac{1}{\alpha+1}\frac{2}{p}>0,
\end{align*}
it follows from Corollary \ref{upper-bound-phi} (iv) 
that
\begin{align*}
 \phi(w)\leq \frac{\alpha}{b+1+\alpha b}=M_1.
\end{align*}
Additionally, by  
Lemma \ref{lower-bound-phi} (iv), 
\begin{align*}
 w\frac{\phi'(w)}{\phi(w)} \geq \frac{(p/2+a+2)b}{b+1}=-M_2.
\end{align*}
By Corollary \ref{cor-minimax} (ii), it follows that the generalized Bayes
estimator is minimax provided
\begin{align}\label{final}
\frac{(n+2)\{ M_1-c(p,n)\}}{1+M_1}+4M_2 \leq 0
\end{align}
but a straightforward calculation shows that
this condition is equivalent to 
\begin{align}
b \geq -(n+2)\frac{c(p,n)-\alpha}{4(p+a+1)(\alpha+1)}.
\end{align}
Hence the generalized Bayes estimator is minimax since 
\eqref{main-b} guarantees \eqref{final} is satisfied. This completes the
 proof. 
\end{proof}
Note: A recent paper by Wells and Zhou (2008) also derives generalized
Bayes estimators in the normal case, some of which have non-monotone
shrinkage functions $\phi$ (in our notation).
The generalized Bayes minimax estimators of Theorem \ref{main-theorem}
(ii) all have non-monotone shrinkage functions and their minimaxity 
cannot be shown by the methods in Wells and Zhou (2008) which rely on
integration by parts. As noted earlier, the assumption that $b<0$
causes the usual integration by parts technique to fail.

\section{Generalized Bayes estimators for spherically symmetric 
distributions}
In this section, we consider generalized Bayes minimax
estimators for spherically symmetric distributions.
As shown in Maruyama (2003) and Maruyama and Strawderman (2005),
the special choice $ b=-a-2$ in the prior given by \eqref{prior}
leads to the separated joint density of $ \theta$ and $ \eta $,
$ \|\theta\|^{-2(p/2+a+1)}\eta^{-a-1+e} $.
This follows since
\begin{align}
&\int_0^1  \exp\left(-\frac{\eta \lambda \|\theta\|^2}{2(1-\lambda)}
 \right)^{p/2}\left( \frac{\eta\lambda}{1-\lambda}
 \right)^{p/2} \lambda^a(1-\lambda)^{b}d\lambda 
\nonumber \\
&= \eta^{-a-1} \int_0^{\infty}t^{p/2+a}\exp \left( -\frac{t}{2} \|\theta\|^2
\right) dt \nonumber \\
& 
\propto \|\theta\|^{-2(p/2+a+1)}\eta^{-a-1}, \label{sep} 
\end{align}
if $p/2+a+1>0$.
Because of this simplification, we
make the assumption $b=-a-2$ throughout the rest of this section.
Under  quadratic loss 
$ \eta (d- \theta)'(d-\theta) $, 
even in the spherically symmetric situation,
the generalized Bayes estimator 
is given by $ E(\eta\theta|X,S)/E(\eta|X,S) $
 and hence we have the generalized Bayes estimator with respect to our prior,
 \begin{align*}
& \frac{\int_{R^{p}} \int_{0}^{\infty} \theta \eta^{(n+p)/2-a+e}
f(\eta\{\|X-\theta\|^2+S\}) \|\theta\|^{-2(p/2+a+1)}
d\eta d\theta}{\int_{R^{p}} \int_{0}^{\infty}
   \eta^{(n+p)/2-a+e} 
f(\eta\{\|X-\theta\|^2+S\}) \|\theta\|^{-2(p/2+a+1)} d\eta d\theta} \\
&=\frac{\int_{R^{p}}  \theta 
(\{\|X-\theta\|^2+S)^{-(n+p)/2+a-e-1} 
\|\theta\|^{-2(p/2+a+1)}
 d\theta \int_0^{\infty} \eta^{(n+p)/2-a+e}f(\eta)d \eta}
{\int_{R^{p}} 
(\{\|X-\theta\|^2+S)^{-(n+p)/2+a-e-1} 
\|\theta\|^{-2(p/2+a+1)} 
 d\theta \int_0^{\infty} \eta^{(n+p)/2-a+e}f(\eta)d \eta} \\
&=\frac{\int_{R^{p}}  \theta 
(\{\|X-\theta\|^2+S)^{-(n+p)/2+a-e-1} 
\|\theta\|^{-2(p/2+a+1)}
 d\theta}
{\int_{R^{p}} 
(\{\|X-\theta\|^2+S)^{-(n+p)/2+a-e-1} 
\|\theta\|^{-2(p/2+a+1)} 
 d\theta }
\end{align*}
if 
\begin{align} \label{moment-condition}
\int_0^{\infty} \eta^{(n+p)/2-a+e}f(\eta)d\eta < \infty.
\end{align}
Note that this does not depend on $f$ and hence is equal to 
the generalized Bayes estimator in the normal case.
In the normal case, as seen in Section 3,
the estimator is well-defined if
$ a> -p/2-1$, $ b>-1$ and $ e> -p/2-n/2-2$.
Since $ a=-b-2$, the inequality $ -p/2-1<a<-1$ is  also satisfied.

We note that Theorem \ref{main-theorem}, together with the general
results of Section 2, imply that the normal theory generalized Bayes
estimators in Theorem \ref{main-theorem} 
remain minimax (but not necessarily generalized
Bayes)
for the entire class of spherically symmetric distributions. The
additional restriction that $b=-a-2$ implies, as noted above, that this
subclass is also generalized Bayes for the entire class of spherically
symmetric distributions. Hence we have
the following result on the minimaxity of the generalized
Bayes estimator with respect to $ \|\theta\|^{-2(p/2+a+1)}\eta^{-a-1+e}$
for the general spherically symmetric case.
\begin{thm}
\begin{enumerate}[(i)]
\item \ Suppose $-p/2-n/2-1<e \leq -n/4-3/2$. When
\begin{align} \label{mg-condition-a}
 -p/2-1<a \leq \frac{c(p,n)(n/2+e)-p/2-1}{1+c(p,n)}
\end{align}
where $c(p,n)=2(p-2)/(n+2)$, the generalized Bayes estimator
is minimax under the moment condition given by \eqref{moment-condition}.
\item \ Suppose $e > \-n/4-3/2$ and $n\geq 2$. When
\begin{align} \label{mg-condition-a-2}
-p/2-1<a \leq a_*
\end{align}
where $a_*$ is the larger solution of the equation $ g(a)=0$ where
\begin{align*}
g(a) &=(2p+2n+4e+4)a^2 \\ & \ + a\{ 2p^2+2pn+12p+7n+4(p+3)e+10\} \\
& \ 4p^2+\{7/2\}pn+6(p+2)e+7n+13p+10,
\end{align*}
the generalized Bayes estimator 
is minimax under the moment condition given by \eqref{moment-condition}.
\end{enumerate}
\end{thm}
Note: Since $ a_*$ will be shown to between $(-2,-1)$ in the proof,
minimaxity of the generalized Bayes estimator with $ -2<a<a_* $ 
($ -1< -a_*-2 < b< 0$), 
which has non-monotone shrinkage factor $\phi$, is new
compared to Maruyama and Strawderman (2005).
\begin{proof} \
First consider that
the upper bound of $a$ for minimaxity in Theorem \ref{main-theorem},
\begin{align*}
u_*(p,n,e)= \frac{c(p,n)(n/2+e)-p/2-1}{1+c(p,n)} .
\end{align*}
A simple calculation gives 
\begin{align*}
u_*(p,n,e) +2 =\frac{(p-2)(n+4e+6)}{2(n+2)\{1+c(p,n)\}}.
\end{align*}
Hence $ u_*(p,n,e)$ is greater than $-2$
if and only if $e>-n/4-3/2$.
To show (i), note that 
when $ -p/2-n/2-1<e \leq -n/4-3/2$,  $ u_*(p,n,e)<-2$ and
\eqref{mg-condition-a} satisfies the sufficient
 condition on $a$ for minimaxity of Theorem \ref{main-theorem} (i). 
Also since $-a-2 \, (=b)$ is nonnegative, 
the  condition on $b$ is also satisfied. Hence part (i) follows.

Next we show (ii).
When $ e > -n/4-3/2$, the generalized Bayes estimator with 
$-p/2-1 <a \leq -2 $ 
is minimax because $ -2< u_*(p,n,e)$  and $-a-2 \, (=b) \geq 0$ guarantees
minimaxity by Theorem \ref{main-theorem} (i).

Finally we consider the case where $ -2<a<-1$ and $ e > -n/4-3/2$.
The condition $ b= -a-2 < \min(n/2+e-a-1, 0)$ in Theorem  
\ref{main-theorem} (ii) is clearly satisfied under the assumption
$n \geq 2$. 
Also a straightforward calculation shows that the inequality
\begin{align} \label{fukuzatsu}
-(n+2)\frac{c(p,n)-\alpha}{4(p+a+1)(\alpha+1)} \leq -a-2 =b
\end{align}
is equivalent to $g(a) \leq 0$.
Note, \eqref{fukuzatsu}, for $a=-1$, is not satisfied because
\begin{align*}
 -a-2 +(n+2)\frac{c(p,n)-\alpha}{4(p+a+1)(\alpha+1)} 
&= -1+(n+2)\frac{c(p,n)-\alpha}{4p(\alpha+1)} \\
& \leq -1+\frac{2(p-2)}{4p(\alpha+1)} \\
& < 0.
\end{align*}
Also  \eqref{fukuzatsu}, for $ a=-2$, is satisfied because
\begin{align*}
 -a-2 +(n+2)\frac{c(p,n)-\alpha}{4(p+a+1)(\alpha+1)} 
&= (n+2)\frac{c(p,n)-\alpha}{4(p-1)(\alpha+1)}  > 0.
\end{align*}
Hence $ g(-2)<0$ and $g(-1)>0$, which guarantees that
the larger solution $a_*$ of the quadratic 
equation $g(a)=0$ should be between
 $-2$ and $-1$. 
Additionally $g(u_{*}(p,n,e))>0$ 
because 
\begin{align*}
 -u_*(p,n,e)-2 +(n+2)\frac{c(p,n)-\alpha}{4(p+u_*(p,n,e)+1)(\alpha+1)} 
&=  -u_*(p,n,e)-2< 0,
\end{align*}
which means that $ a_* $ is smaller than $u_*(p,n,e)$.
Hence $ -2<a \leq a_*$ leads to  minimaxity. This completes the proof of
 (ii). 
\end{proof}

\end{document}